\documentclass[12pt]{amsart}
\usepackage{txfonts}
\usepackage{mathrsfs}
\usepackage{amsmath}
\usepackage{amssymb, amsmath}
\usepackage{graphicx}

\newcommand{\newcom}{\newcommand}

\newcom{\al}{\alpha}
\newcom{\be}{\beta}
\newcom{\eps}{\epsilon}
\newcom{\veps}{\varepsilon}
\newcom{\ga}{\gamma}
\newcom{\Ga}{\Gamma}
\newcom{\ka}{\kappa}
\newcom{\Lam}{\Lambda}
\newcom{\lam}{\lambda}
\newcom{\Om}{\Omega}
\newcom{\om}{\omega}
\newcom{\Si}{\Sigma}
\newcom{\si}{\sigma}
\newcom{\tht}{\theta}
\newcom{\dtri}{\nabla}
\newcom{\tri}{\triangle}
\newcom{\oo}{\infty}
\newcom{\vphi}{\varphi}
\newcom{\cB}{{\mathcal B}}
\newcom{\cC}{{\mathcal C}}
\newcom{\cD}{{\mathcal D}}
\newcom{\cF}{{\mathcal F}}
\newcom{\cL}{{\mathcal L}}
\newcom{\cM}{{\mathcal M}}
\newcom{\cP}{{\mathcal P}}
\newcom{\cS}{{\mathcal S}}
\newcom{\cQ}{{\mathcal Q}}
\newcom{\cT}{{\mathcal T}}
\newcom{\cY}{{\mathcal Y}}
\newcom{\cZ}{{\mathcal Z}}
\newcom{\R}{\mathbb R}
\newcom{\T}{\mathbb T}
\newcom{\N}{\mathbb N}
\newcom{\Z}{\mathbb Z}
\newcom{\C}{\mathbb C}
\newcom{\E}{\mathbb E}

\newcom{\f}{\frac}
\newcom{\di}{\displaystyle\int}
\newcom{\ds}{\displaystyle\sum}
\newcom{\dl}{\displaystyle\lim}
\newcom{\ov}{\overline}
\newcom{\sset}{\subset}
\newcom{\wt}{\widetilde}
\newcom{\pa}{\partial}
\newcom{\p}{\partial}
\newcom\na{\nabla}
\newcom{\suml}{\sum\limits}
\newcom{\supl}{\sup\limits}
\newcom{\intl}{\int\limits}
\newcom{\infl}{\inf\limits}
\newcom{\disp}{\displaystyle}
\newcom{\non}{\nonumber}
\newcom{\no}{\noindent}
\newcom{\QED}{$\square$}

\def\div{\mathop{\rm div}\nolimits}

\def\ef{\hphantom{MM}\hfill\llap{$\square$}\goodbreak}

\def\eqdefa{\buildrel\hbox{\footnotesize def}\over =}

\newtheorem{athm}{\bf \t}[section]
\newenvironment{thm} [1] {\def\t{#1}\begin{athm} \bf \rm} {\end{athm}}
\newcom{\bthm}{\begin{thm}}\newcom{\ethm}{\end{thm}}

\newtheorem{proposition}{Proposition}[section]

\newcom{\beq}{\begin{equation}}
\newcom{\eeq}{\end{equation}}
\newcom{\ben}{\begin{eqnarray}}
\newcom{\een}{\end{eqnarray}}
\newcom{\beno}{\begin{eqnarray*}}
\newcom{\eeno}{\end{eqnarray*}}
\newcom{\bali}{\begin{aligned}}
\newcom{\eali}{\end{aligned}}

\numberwithin{equation}{section}

\begin{document}

\title[compressible Navier-Stokes equations]
{A Beale-Kato-Majda Blow-up criterion for the 3-D compressible
Navier-Stokes equations}

\author{Yongzhong Sun}
\address{Department of Mathematics, Nanjing University, 210093, P. R. China}
\email{sunyongzhong@163.com}

\author{Chao Wang}
\address{Academy of Mathematics $\&$ Systems Science, CAS, Beijing 100190, P. R. China}
\email{wangchao@amss.ac.cn}

\author{Zhifei Zhang}
\address{School of Mathematical Sciences, Peking University, 100871, P. R. China}
\email{zfzhang@math.pku.edu.cn}


\date{January 6, 2010}


\begin{abstract}
We prove a blow-up criterion in terms of  the upper bound of the
density for the strong solution to the 3-D compressible Navier-Stokes equations.
The initial vacuum is allowed. The main ingredient of the proof is \textit{a priori} estimate
for an important quantity under the assumption that the density is upper
bounded, whose divergence can be viewed as the effective viscous flux.
\end{abstract}

\maketitle

\section{Introduction}

In this paper, we consider the isentropic compressible Navier-Stokes
system in three dimensional space. The system reads
\begin{eqnarray}\label{NS}
\left\{\begin{aligned}
&\partial_t \rho+\div(\rho u)=0,\hspace{4pt} \text{in}\hspace{4pt} (0,T)\times\Om,\\
&\partial_t (\rho u)+\div(\rho u\otimes u)-Lu+\nabla
p=0,\hspace{4pt} \text{in}\hspace{4pt} (0,T)\times\Om,
\end{aligned}\right.
\end{eqnarray}
together with the initial-boundary conditions
\begin{eqnarray}\label{initial}
&&(\rho(t,x),u(t,x))|_{t=0}=(\rho_0(x),u_0(x)),\quad \text{in}\hspace{4pt} \Om,\\
&&u(t,x)=0,\quad \text{on}\hspace{4pt} (0,T)\times\p\Om.
\end{eqnarray}
Here $\Om$ is either $\mathbb{R}^3$ or a bounded domain in $\R^3$,
$\rho$ and $u$ are the density and velocity of the fluid
respectively, $p=a\rho^{\gamma}$ with $\gamma>1$ is the pressure.
The Lam\'{e} operator $L$ is defined by
\begin{eqnarray*}
Lu=\mu\Delta u+(\mu+\lambda)\nabla\div u,
\end{eqnarray*}
with constant viscosity coefficients $\mu$ and $\lambda$ satisfying
\begin{eqnarray}\label{Co}
\mu>0,\hspace{4pt} 3\lambda+2\mu\geq 0.
\end{eqnarray}

In the absence of vacuum for the initial density, the local
existence of strong solution as well as the global existence of
strong solution and weak solution with the initial data close to an
equilibrium state were well developed, see \cite{Mat,Nash,Sol,Hoff,Dan,Chen}
and references therein. The global existence of weak solution for
large initial data was first solved by P. L. Lions in \cite{Lions}
for $\ga\ge \f95$. E. Feireisl, A. Novotn\'{y} and H.
Petzeltov\'{a} \cite{Fei} extended Lions's result to the case of
$\ga>\f32$. S. Jiang and P. Zhang \cite{Jiang1,Jiang2} proved the global existence
of weak solution for any $\ga>1$ for the spherically
symmetric or axisymmetric initial data. However, the regularity and uniqueness of
weak solutions are completely open even in the case of two
dimensional space. The only known result is the work of Kazhikhov
and Va\u{\i}gant \cite{Vaigant}, where they proved the global existence of
strong solution for the system (\ref{NS}) in $\Om=\mathbb{T}^2$
under the assumption that $\mu$ is a constant and
$\lambda=\rho^{\beta}$ with $\beta>3$. On the other hand, when the
initial density is compactly supported, Z. Xin \cite{Xin} proved that
smooth solution will blow up in finite time in the whole space.

To proceed we introduce some notations for the standard homogeneous
and inhomogeneous Sobolev spaces.
\[
D^{k,r}(\Om)\eqdefa\{u\in L^1_{loc}(\Om):\|\nabla^k
u\|_{L^r(\Om)}<\infty\},
\]
\[
W^{k,r}(\Om)\eqdefa L^r(\Om)\cap
D^{k,r}(\Om),H^k(\Om)=W^{k,2}(\Om),D^k(\Om)=D^{k,2}(\Om),
\]
\[
D^1_0(\Om)\eqdefa \{u\in L^6(\Om):\|\nabla u\|_{L^2(\Om)}<\infty\,\,
and\,\, u=0 \,\, on\,\, \pa\Om\},
\]
\[
H^1_0(\Om)\eqdefa L^2(\Om)\cap
D^1_0(\Om),\|u\|_{D^{k,r}(\Om)}=\|\nabla^k u\|_{L^r(\Om)}.
\]

When the initial vaccuum is allowed, the local well-posedness and
blow-up criterion for strong solutions to the compressible
Navier-Stokes equations  were established in a series of
papers \cite{CCK,CK1,CK3} by  Cho, Choe and Kim. Here we write down
one of those results.
\bthm{Theorem}\label{thmCCSS} Let $\Omega$ be
a bounded smooth domain or $\mathbb{R}^3$ and $q\in (3,6]$. Suppose
that $\rho_0\geq 0$ and belongs to $W^{1,q}(\Omega)\cap H^1(\Om)\cap
L^1(\Om)$, $u_0\in D^1_0(\Omega)\cap D^2(\Omega)$ with the following
compatibility condition satisfied,
\begin{eqnarray}\label{CompC}
-\mu\Delta u_0-(\mu+\lambda)\nabla\div u_0+\nabla
p(\rho_0)=\sqrt{\rho_0}g,
\end{eqnarray}
for some vector field $g\in L^2(\Omega)$. Then there exist a time
$T\in (0,\infty]$ and a unique strong solution $(\rho,u)$ to
(\ref{NS}) such that
\begin{eqnarray*}
\rho\in C([0,T),H^1\cap W^{1,q}(\Omega)),u\in
C([0,T),D^2(\Omega))\cap L^2(0,T;D^{2,q}(\Omega)).
\end{eqnarray*}
Moreover, let $T^*$ be a maximal existence time of the solution. If
$T^*<\infty$, then there holds
\begin{eqnarray}\label{blowup1}
\limsup\limits_{t\uparrow
T^*}\left(\|\rho(t)\|_{W^{1,q}(\Omega)}+\|u(t)\|_{D_0^1(\Omega)}\right)=\infty.
\end{eqnarray}
\ethm

Since the initial vacuum is allowed,
it is then important to investigate the possible blow-up mechanism of the solution.
In their recent works \cite{Huang1,Huang2}, X. Huang and Z. Xin established
a Beale-Kato-Majda blow up criterion for the above strong solution. More
precisely,
\bthm{Theorem}\label{HX} Assume that the coefficients of the
operator $L$ satisfies (\ref{Co}) and moreover,
\begin{eqnarray}\label{coefficient}
\lambda<7\mu.
\end{eqnarray}
Let $(\rho,u)$ be the strong solution constructed in Theorem \ref{thmCCSS}
and $T^*$ be a maximal existence time. If $T^*<\infty$, then
\begin{eqnarray}\label{blowup2}
\lim\limits_{T\uparrow T_*}\|\nabla
u\|_{L^{1}(0,T,L^{\infty}(\Omega))}=\infty.
\end{eqnarray}
\ethm
Recently, J. Fan, S. Jiang and Y. Ou \cite{Fan} also obtained a similar result for the compressible heat-conductive flows.
On the other hand, for the 2D compressible Navier-Stokes equations in $\T^2$, B. Desjardins \cite{Des} proved
more regularity of weak solution under the assumption that the density is upper bounded; Very recently,
L. Jiang and Y. Wang \cite{Jiangw}, Y. Sun and Z. Zhang \cite{Sun} obtained a blow-up criterion
in terms of the upper bound of the  density for the strong solution. In \cite{Sun},
the initial vacuum is allowed and the domain includes the bounded domain.
Note that the $L^{1}(0,T,L^{\infty}(\Omega))$ bound for $\nabla u$
immediately implies the upper bound for the density $\rho$.

The purpose of this paper is to obtain a Beale-Kato-Majda blow-up criterion in terms of the upper bound of the
density for the 3-D compressible Navier-Stokes equations. Our main result is stated as follows.

\bthm{Theorem}\label{main}
Assume that $(\rho,u)$ is the strong solution constructed in Theorem
\ref{thmCCSS}. Let $\mu,\lambda$ be as in Theorem \ref{HX} and $T^*$ be a maximal existence time of the solution.
If $T^*<\infty$, then we have
\begin{eqnarray}\label{blowup}
\limsup\limits_{T\uparrow
T^*}\|\rho(t)\|_{L^{\infty}(0,T;L^\infty(\Omega))}=\infty.
\end{eqnarray}
\ethm

\bthm{Remark}
This result seems surprise, if we compare with the incompressible Navier-Stokes equations where the density is
a constant. It is well-known that if we have some kind of control for the pressure, the Leray weak solution is in fact
smooth for the incompressible Navier-Stokes equations, see \cite{Ber}. For the compressible Navier-Stokes equations,
the pressure is determined by the density, the bound of the density
thus implies a bound for the pressure. From this viewpoint, our result seems natural.
\ethm

\bthm{Remark}
In a forthcoming paper, we will extend similar result
to the compressible heat-conductive flows.
\ethm

Let us conclude this section by introducing the main idea of our proof. First of all, if the density is upper bounded,
we can obtain a high integrability of the velocity, see Lemma \ref{uLr}. This bound can be used to
control the nonlinear term. The trouble is to control the density, which satisfies a transport equation.
In order to propagate the regularity of the density, it is necessary to require that the velocity is bounded
in $L^1(0,T;W^{1,\infty}(\Om))$. On the other hand, we have to obtain a priori bound of $\na \rho$
in order to prove $u\in L^1(0,T;W^{1,\infty}(\Om))$. To overcome this difficulty, we
introduce an important quantity $w$ defined by $w=u-v$,
where $v$ is the solution of  Lam\'{e} system
\beno \left\{\begin{aligned}
&\mu\Delta v+(\lambda+\mu)\nabla \div v=\nabla p(\rho)\hspace{4pt}\text{in}\hspace{4pt} \Omega,\\
&v(x)=0\hspace{4pt}\text{on} \hspace{4pt}\partial\Omega.
\end{aligned}\right.
\eeno
In the case of $\Om=\R^3$, $(\lambda+2\mu)\textrm{div} w=(\lambda+2\mu)\textrm{div} u-p
\eqdefa G$. It is well known that $G$ is called the effective viscous flux, which plays
an important role in the existence theory of weak solution. A key point is that  we can obtain the
better regularity of $w$ than $u$ under the only assumption that the density is upper bounded.
More precisely, we proved that $\na^2 w\in L^2(0,T;L^6(\Om))$, which combined with
the bound of the density implies that $\na u\in L^2(0,T;L^\infty(\Om)+L^\infty(0,T;BMO(\Om))$.
This bound still does not imply that $\na u$ is bounded in $L^1(0,T;L^\infty(\Om))$. We need to
introduce the second key ingredient: a logarithmic estimate for
Lam\'{e} system. Then the result can be deduced by combining the
above two estimates into the energy estimates for the density.

\section{Preliminaries}

Consider the following boundary value problem for the Lam\'{e}
operator $L$
\begin{eqnarray}\label{eq:Lame}
\left\{\begin{aligned}
&\mu\Delta U+(\mu+\lambda)\nabla\div U=F, \hspace{4pt}\text{in}\hspace{4pt} \Omega,\\
&U(x)=0,\hspace{4pt}\text{on} \hspace{4pt}\partial\Omega.
\end{aligned}\right.
\end{eqnarray}
Here $U=(U_1,U_2,U_3)$, $F=(F_1,F_2,F_3)$. It is well known that
under the assumption (\ref{Co}), (\ref{eq:Lame}) is a strongly
elliptic system. If $F\in W^{-1,2}(\Omega)$, then there exists an
unique weak solution $U\in D_0^1(\Omega)$. We begin with recalling
various estimates for this system in $L^q(\Omega)$ spaces.

\bthm{Proposition}\label{proL} Let $q\in (1,\infty)$ and $U$ be a
solution of (\ref{eq:Lame}). There exists a constant $C$ depending
only on $\lambda,\mu,q$ and $\Omega$ such that the following estimates hold.

(1)\,If $F\in L^q(\Om)$, then
\begin{eqnarray}\label{E1}
\left\{\begin{aligned}
&\|D^2U\|_{L^q(\R^3)}\leq C\| F\|_{L^q(\R^3)},\\
&\|U\|_{W^{2,q}(\Omega)}\leq C\| F\|_{L^q(\Omega)};\quad \text{if}\hspace{2pt} \Om\hspace{2pt} \text{is a bounded domain}.
\end{aligned}\right.
\end{eqnarray}

(2)\,If $F\in W^{-1,q}(\Omega)$(i.e., $F=\div f$ with
$f=(f_{ij})_{3\times 3}, f_{ij}\in L^{q}(\Omega)$), then
\begin{eqnarray}\label{E2}
\left\{\begin{aligned}
&\|DU\|_{L^q(\R^3)}\leq C\| f \|_{L^q(\R^3)},\\
&\|U\|_{W^{1,q}(\Omega)}\leq C\| f\|_{L^q(\Omega)};\quad \text{if}\hspace{2pt} \Om \hspace{2pt}\text{is a bounded domain}.
\end{aligned}\right.
\end{eqnarray}

(3)\,If $F=\div f$ with $f_{ij}=\partial_k h^k_{ij}$ and
$h^k_{ij}\in W_0^{1,q}(\Omega)$ for $i,j,k=1,2,3$, then
\begin{eqnarray}
\label{E3} \| U\|_{L^{q}(\Omega)} \leq C\| h \|_{L^q(\Omega)}.
\end{eqnarray}

\ethm

\noindent{\bf Proof.}\, In the case when $\Om$ is a bounded domain,
the estimates (\ref{E1}) and (\ref{E2}) are classical for strongly
elliptic systems, see for example \cite{ADN}. The estimate
(\ref{E3}) can be proved by a duality argument with the help of
(\ref{E1}). In the case of $\Om=\mathbb{R}^3$, one can give an
explicit representation formula for the solution as follows. Taking
divergence on both sides of (\ref{eq:Lame}), one finds
\begin{eqnarray*}
\div U=\f{1}{\lambda+2\mu}\Delta^{-1}\div F.
\end{eqnarray*}
Substituting this into (\ref{eq:Lame}) gives us
\begin{eqnarray*}
\Delta
U=\f{1}{\mu}F-\f{\lambda+\mu}{\mu(\lambda+2\mu)}\nabla\Delta^{-1}\div
F.
\end{eqnarray*}
Denote the Riesz transform $R=(R_1,R_2,R_3)=\nabla\Delta^{-1/2}$.
Then
\begin{eqnarray*}
\Delta U=\f{1}{\mu}F-\f{\lambda+\mu}{\mu(\lambda+2\mu)}R(R\cdot F).
\end{eqnarray*}
Hence for $i,j,k=1,2,3$,
\begin{eqnarray*}\label{represent1}
\partial_{ij} U_k=\f{1}{\mu}R_iR_jF_k-\f{\lambda+\mu}{\mu(\lambda+2\mu)}R_iR_jR_k(R\cdot F).
\end{eqnarray*}
The classical $L^q(\mathbb{R}^3)$-boundedness for Riesz transform
gives
\begin{eqnarray*}
\|D^2U\|_{L^q(\R^3)}\leq
C(q)\f{2\lambda+3\mu}{\mu(\lambda+2\mu)}\|
F\|_{L^q(\R^3)}.
\end{eqnarray*}
Similar argument gives the estimates (\ref{E2}) and (\ref{E3}).\ef

We need an endpoint estimate for $L$ in the case $q=\infty$. Let
$BMO(\Om)$ stand for the John-Nirenberg's space of bounded mean
oscillation whose norm is defined by
\beno
\|f\|_{BMO(\Omega)}\eqdefa\|f\|_{L^2(\Omega)}+[f]_{BMO}, \eeno
with
\beno
[f]_{BMO(\Omega)}\eqdefa\sup\limits_{x\in {\Om},r\in (0,d)}\fint\limits_{\Om_r(x)}|f(y)-f_{\Om_r(x)}|dy,\\
f_{\Om_r(x)}=\fint\limits_{\Om_r(x)}f(y)dy=\f{1}{|\Om_r(x)|}\int\limits_{\Om_r(x)}f(y)dy.
\eeno Here $\Om_r(x)=B_r(x)\cap\Om$, $B_r(x)$ is the ball with
center $x$ and radius $r$ and $d$ is the diameter of $\Om$.
$|\Om_r(x)|$ denotes the Lebesque measure of $\Om_r(x)$. Note that
$$[f]_{BMO(\Omega)}\leq 2\|f\|_{L^{\infty}(\Omega)}.$$

\bthm{Proposition}\label{propend} If $F=\div f$ with
$f=(f_{ij})_{3\times 3}, f_{ij}\in L^{\infty}(\Omega)\cap L^2(\Om)$,
then $\na U\in BMO(\Omega)$ and there exists a constant $C$ depending
only on $\lambda,\mu$ and $\Omega$ such that
\begin{eqnarray}
\label{E4} \|\nabla U\|_{BMO(\Omega)} \leq C\left(\| f
\|_{L^{\infty}(\Omega)}+\|f\|_{L^2(\Omega)}\right).
\end{eqnarray}
\ethm \noindent{\bf Proof.} When $\Om$ is a bounded domain,the
estimate (\ref{E4}) can be found in \cite{Acq} for a more general
setting. Now if $\Om=\mathbb{R}^3$ we use the representation formula
for $\nabla U$. Since
\begin{equation*}
\Delta U=\f{1}{\mu}\div
f-\f{\lambda+\mu}{\mu(\lambda+2\mu)}\nabla\Delta^{-1}\div\div
f=\f{1}{\mu}\div f-\f{\lambda+\mu}{\mu(\lambda+2\mu)}\nabla G,
\end{equation*}
with $G=\sum^3_{i,j=1}R_iR_jf_{ij}$. For $k,l=1,2,3$,
\begin{equation*}
\partial_{k} U_l=\f{1}{\mu}R_k\sum^3_{j=1}R_jf_{lj}-\f{\lambda+\mu}{\mu(\lambda+2\mu)}R_kR_lG.
\end{equation*}
By the Fefferman-Stein's classical result on $BMO$-boundedness of
singular integral operators \cite{Stein}, there exists an absolute
constant $C>0$ such that
\begin{eqnarray*}
[\nabla U]_{BMO(\R^3)} \leq
C\f{2\lambda+3\mu}{\mu(\lambda+2\mu)}\|
f\|_{L^{\infty}(\Om)}.
\end{eqnarray*}
This inequality combined with  (\ref{E2}) with $q=2$ yields
(\ref{E4}).\ef

In the next lemma, we will give a variant of the Brezis-Waigner's
inequality \cite{Brezis}. To our knowledge, such a kind of
inequality was first established in \cite{Kozono} in the case of
$\Om=\R^3$. For the reader's convenience, we will give a proof in
the case when $\Om$ is a bounded Lipschitz domain, see also
\cite{Sun}.

\bthm{Lemma}\label{embedding2} Let $\Om=\mathbb{R}^3$ or be a
bounded Lipschitz domain and $f\in W^{1,q}(\Om)$ with $q\in
(3,\infty)$. There exists a constant $C$ depending on $q$ and the
Lipshitz property of $\Om$ such that
\begin{eqnarray}\label{embedineq2}
\|f\|_{L^{\infty}(\Omega)}\leq
C\left(1+\|f\|_{BMO(\Omega)}\ln\left(e+\|\nabla
f\|_{L^{q}(\Omega)}\right)\right).
\end{eqnarray}
\ethm

\noindent{\bf Proof.}\,\,First note that for a Lipschitz domain, the
following so-called $A$-property holds:

There exist two constants $A\geq 1$ and $r_0\in(0,d)$ such that for
any $r\in (0,r_0)$ and $x\in \Om$, \beno |\Om_r(x)|\leq |B_r(x)|\leq
A|\Om_r(x)|. \eeno Without loss of generality we assume $r_0\leq 1$.

First we give an estimate for $|f_{\Om_r(x)}|$ with $0<r<r_0$ and
$x\in \Om$. If $r \ge\f{1}{2}r_0$, then \beno
\left|f_{\Om_r(x)}\right|\leq\f{1}{|\Om_r(x)|}\int\limits_{\Om_r(x)}\left|f(y)\right|dy\leq
C\|f\|_{L^2(\Om)}. \eeno
If $r<\f{1}{2}r_0$, then there exists some integer $k\ge 1$ such that
\beno \f{r_0}{2^{k+1}}\leq
r<\f{r_0}{2^{k}},\hspace{4pt} k\leq C(1+|\ln r|).
\eeno
Denoting $\Om_j=\Om_{2^jr}(x)$ for $j=0,1,\cdots,k$, we have
\beno
\left|f_{\Om_r(x)}\right|&\leq&
\sum^{k}_{j=1}\left|f_{\Om_{j-1}}-f_{\Om_j}\right|+\left|f_{\Om_k}\right|\\
&\leq& \sum^{k}_{j=1}\fint_{\Om_{j-1}}\left|f(y)-f_{\Om_j}\right|dy+C\left\|f\right\|_{L^2(\Om)}\\
&\leq& 2^NA\sum^{k}_{j=1}\fint_{\Om_{j}}\left|f(y)-f_{\Om_j}\right|dy+C\left\|f\right\|_{L^2(\Om)}\\
&\leq& C k [f]_{BMO(\Om)}+C\|f\|_{L^2(\Om)}\leq C(1+|\ln
r|)\|f\|_{BMO(\Om)}.
\eeno
We conclude that there exists a constant $C=C(A,r_0,N)$ such that
\beno \left|f_{\Om_r(x)}\right|\leq
C(1+|\ln r|)\|f\|_{BMO(\Om)},
\eeno
which together with Sobolev embedding theorem in a Lipschitz domain \cite{Adam} ensures that
for any fixed $x\in\Om$ and  small enough $\varepsilon>0$ we have
\beno
|f(x)|\leq |f(x)-f_{\Om_{\varepsilon}(x)}|+|f_{\Om_{\varepsilon}(x)}|\leq
C\left(\varepsilon^{1-\f{N}{q}}\|f\|_{W^{1,q}(\Om)}+(1+|\ln\varepsilon|)\|f\|_{BMO(\Om)}\right).
\eeno
A suitable choice of $\varepsilon$ yields the inequality (\ref{embedineq2}).\ef

In the subsequent context we will use $L^{-1}F$ to denote the unique
solution $U$ of the Lam\'{e} system (\ref{eq:Lame}).

\section{A priori estimates for the effective viscous flux}

In what follows, we assume that $(\rho,u)$ is a strong solution of
(\ref{NS}) in $[0,T)$ with the regularity stated in Theorem
\ref{thmCCSS}.

Standard energy estimates yields that for any $t\in [0,T)$,
\begin{eqnarray*}
&&\|\rho(t)\|_{L^{1}(\Omega)}\leq \|\rho_0\|_{L^{1}(\Omega)},\\
&&\|\rho(t)\|^\gamma_{L^{\gamma}(\Omega)}+\|\rho
|u|^2(t)\|_{L^1(\Omega)}+\|\nabla u\|^2_{L^2((0,t)\times\Om)}\\
&&\leq
C\big( \|\rho_0\|^\gamma_{L^{\gamma}(\Omega)}+\|\rho_0
|u_0|^2\|_{L^1(\Omega)}\big).
\end{eqnarray*}
Note that by the assumption on $\rho_0, u_0$,
\begin{eqnarray*}
\|\rho_0\|^{\gamma}_{L^{\gamma}(\Omega)}\leq
\|\rho_0\|^{\gamma-1}_{L^{\infty}(\Omega)}\|\rho_0\|_{L^{1
}(\Omega)},\quad \|\rho_0 |u_0|^2\|_{L^1(\Omega)}\leq
\|\rho_0\|_{L^{3/2}(\Omega)}\| u_0\|^2_{L^6(\Omega)}.
\end{eqnarray*}
We thus have the following bounds
\begin{eqnarray}\label{eq:energy}
\|\rho\|_{L^{\infty}(0,T;L^{1}(\Omega))},\quad
\|\sqrt{\rho}u\|_{L^{\infty}(0,T;L^2(\Omega))},\quad \|\nabla
u\|_{L^{2}(0,T;L^2(\Omega))}\leq C.
\end{eqnarray}
Here $C$ depends only on $\mu,\lambda,\gamma,a$ and $\rho_0, u_0$.

In what follows the dependence of the constant $C$ on
$\mu,\lambda,\gamma,a$ and $\Om$ will not be mentioned.

The following lemma is the first key step, whose argument comes from
\cite{Hoff} and \cite{Huang2}.
\bthm{Lemma}\label{uLr} Assume that
$\mu<7\lambda$ and the density $\rho$ satisfies
\ben\label{M}
\|\rho\|_{L^\infty(0,T;L^\infty(\Om))}\le M. \een
There exists $r\in (3,6)$ such that $\rho |u|^r\in L^{\infty}(0,T;L^1(\Omega))$ with
$$
\|\rho|u|^r\|_{L^{\infty}(0,T;L^1(\Omega))}\leq C.
$$
Here $C$ depends on $T, \|\rho_0\|_{L^\infty(\Om)}, \|\nabla
u_0\|_{L^2(\Om)},M$. \ethm

\noindent{\bf Proof.} Multiplying the second equation of (\ref{NS})
by $r|u|^{r-2}u$, and integrating the resulting equation on $\Omega$
to obtain
\begin{eqnarray}\label{eq1}
&&\frac{d}{dt}\int\limits_{\Omega}\rho |u|^rdx+\int\limits_{\Omega}r|u|^{r-2}\left(\mu|\nabla u|^2+(\lambda+\mu)(\div u)^2\right)\nonumber\\
&&\qquad+r(r-2)\left(\mu|u|^{r-2}|\nabla|u||^2+(\lambda+\mu)(\div u)|u|^{r-3}u\cdot\nabla|u|\right)dx\nonumber\\
&&\quad=\int\limits_{\Omega}r p(\rho)\div(|u|^{r-2}u)dx
\end{eqnarray}

By using the fact $|\nabla u|\geq |\nabla|u||$, the term in the
second integrand can be estimated from below by
\begin{eqnarray*}
&&r|u|^{r-2}\Big[\mu|\nabla u|^2+(\lambda+\mu)(\div u)^2+(r-2)\mu|\nabla|u||^2\\
&&\quad-(\lambda+\mu)(r-2)|\nabla|u|||\div u|\Big]\\
&&\ge r|u|^{r-2}\big[\mu|\nabla u|^2+(\lambda+\mu)\big(\div
u-\f{r-2}{2}|\nabla|u||\big)^2-(\lambda+\mu)\f{(r-2)^2}{4}|\nabla|u||^2\\
&&\qquad\qquad+(r-2)\mu
|\nabla|u||^2 \big]\\
&&\ge r|u|^{r-2}\big[\mu|\nabla
u|^2+(r-2)\big(\mu-(\lambda+\mu)\f{r-2}{4}\big)|\nabla|u||^2 \big]
\end{eqnarray*}
Recalling that $\lambda<7\mu$, there exists $r\in (3,6)$ such that
the last term is greater than
\begin{eqnarray*}
c |u|^{r-2}|\nabla u|^2.
\end{eqnarray*}
On the other hand, because of $\|\rho\|_{L^\infty}\leq M$, we find that the right-hand side of (\ref{eq1})
is controlled by
\begin{eqnarray*}
C\int\limits_{\Omega}\rho^{\frac{r-2}{2r}}|u|^{r-2}|\nabla u|dx \leq
\epsilon\int\limits_{\Omega}|u|^{r-2}|\nabla
u|^2dx+{\frac{C}{\epsilon}}\left(\int\limits_{\Omega}\rho|u|^{r}dx\right)^{\frac{r-2}{r}}.
\end{eqnarray*}
Taking $\epsilon=\frac{c}{2}$ to yield that
\begin{eqnarray*}
\frac{d}{dt}\int\limits_{\Omega}\rho |u|^rdx\leq
C\left(\int\limits_{\Omega}\rho|u|^{r}dx\right)^{\frac{r-2}{r}},
\end{eqnarray*}
which together with the following bound
$$\|\rho_0 |u_0|^r\|_{L^1(\Om)}\leq
\|\rho_0\|_{L^{\f{6}{6-r}}(\Om)}\|u_0\|^r_{L^6(\Om)}\leq
C\|\rho_0\|_{L^{\f{6}{6-r}}(\Om)}\|\nabla u_0\|^r_{L^2(\Om)},$$
implies the desired estimate.\ef
\vspace{0.1cm}

Now for each $t\in [0,T)$, we denote $v(t,x)\eqdefa L^{-1}\na p(\rho)$. That
is, $v(t)$ is the solution of \ben\label{eqv} \left\{\begin{aligned}
&\mu\Delta v+(\lambda+\mu)\nabla \div v=\nabla p(\rho)\hspace{4pt}\text{in}\hspace{4pt} \Omega,\\
&v(t,x)=0\hspace{4pt}\text{on} \hspace{4pt}\partial\Omega.
\end{aligned}\right.
\een
Thanks to Proposition \ref{proL},  for any $q\in (1,\infty)$, there exists
a constant $C$ independent of $t$ such that
\ben\label{estv}
\begin{array}{ll}
&\|\nabla v(t)\|_{L^q(\Omega)}\leq C\|p(\rho(t))\|_{L^q(\Omega)},\\
&\|\nabla^2 v(t)\|_{L^q(\Omega)}\leq C\|\nabla
p(\rho(t))\|_{L^q(\Omega)}.
\end{array}
\een

Now let us introduce an important quantity
$$
w=u-v,
$$
whose divergence can be viewed as the effective viscous flux.

An important observation is that this quantity possesses more regularity information than $u$
does under the assumption that the density is upper bounded. More precisely,

\bthm{Proposition}\label{lemmaw1} Under the assumption (\ref{M}), we
have
\ben\label{estw1} \|\nabla
w\|_{L^{\infty}(0,T;L^2(\Om))},\hspace{4pt} \|\rho^{\f12}\partial_t
w\|_{L^2((0,T)\times\Om)},\hspace{4pt} \|\nabla^2
w\|_{L^2((0,T)\times\Om)}\leq C. \een Here the constant $C$ depends
on $\|\rho_0\|_{L^\infty(\Om)}, \|\nabla u_0\|_{L^2(\Om)}, M, T$.
\ethm

\noindent{\bf Proof.}\, By using the continuity equation, we find
that $w$ satisfies
\ben\label{eqw} \left\{\begin{aligned}
&\rho\partial_t w-\mu\Delta w-(\lambda+\mu)\nabla \div w=\rho F,\hspace{4pt}\text{in}\hspace{4pt} (0,T)\times\Om,\\
&w(t,x)=0\hspace{4pt}\text{on}
\hspace{4pt}[0,T)\times\partial\Omega, \hspace{4pt}
w(0,x)=w_0(x),\hspace{4pt} \text{in}\hspace{4pt} \Omega,
\end{aligned}\right.
\een with $w_0(x)=u_0(x)+v_0(x)$ and \beno
F&=&-u\cdot\nabla u-L^{-1}\nabla(\pa_tp(\rho))\non\\
&=&-u\cdot\nabla u+L^{-1}\nabla\div[p(\rho)u]+L^{-1}\nabla[(\rho
p'(\rho)-p(\rho))\div u]. \eeno

Multiplying the first equation of (\ref{eqw}) by $\partial_t w$ and integrating the
resulting equation over $\Omega$ to obtain ,
\beno
\frac{d}{dt}\int\limits_{\Omega}\mu|\nabla w|^2+(\lambda+\mu)|\div
w|^2dx+\int\limits_{\Omega}\rho|\partial_t w|^2dx
=\int\limits_{\Om}\rho F\cdot\pa_t w dx,
\eeno
which together with H\"{o}lder inequality and Young's inequality gives
\ben\label{eq:w-energy}
&&\frac{d}{dt}\int\limits_{\Omega}\mu|\nabla
w|^2+(\lambda+\mu)|\div w|^2dx+\f{1}{2}\int\limits_{\Omega}\rho|\partial_t w|^2dx\non\\
&&\leq \f{1}{2}\|\sqrt{\rho}F\|^2_{L^2(\Om)}.
\een

Now let us estimate $\|\sqrt{\rho}F\|^2_{L^2(\Om)}$. We get by Lemma
\ref{uLr} and (\ref{estv}) that
\beno
\|\sqrt{\rho} u\cdot\nabla
u\|_{L^2(\Om)}
&\leq& C\|\rho^{\f1r} u\|_{L^r(\Om)}\|\nabla u\|_{L^{\f{2r}{r-2}}(\Om)}\\
&\leq& C\|\rho^{\f1r} u\|_{L^r(\Om)}\big(\|\nabla w\|_{L^{\f{2r}{r-2}}(\Om)}+\|\nabla v\|_{L^{\f{2r}{r-2}}(\Om)}\big)\non\\
&\leq& C_\epsilon\|\nabla w\|_{L^2(\Om)}+\epsilon\|\nabla^2
w\|_{L^2(\Om)}+C. \eeno Here we use the interpolation inequality
\beno \|f\|_{L^q(\Omega)}\leq
C_\epsilon\|f\|_{L^2(\Omega)}+\epsilon\|\nabla
f\|_{L^2(\Omega)},\hspace{4pt} 2\le q<6.
 \eeno
We infer from Proposition \ref{proL} that
\beno &&\|\sqrt{\rho}
L^{-1}\nabla\div[p(\rho)u]\|_{L^2(\Om)}\leq
C\|p(\rho)u\|_{L^2(\Om)}\leq C\|\sqrt{\rho}u\|_{L^2(\Om)}\leq C,\\
&&\|\sqrt{\rho} L^{-1}\nabla (\rho p'-p)\div u\|_{L^2(\Om)}\\
&&\leq \|\sqrt{\rho}\|_{L^3(\Om)}\|L^{-1}\nabla(\rho p'-p)\div u\|_{L^6(\Om)}\\
&&\leq C\|\nabla L^{-1}\nabla(\rho p'-p)\div u\|_{L^2(\Om)}\leq
C\|\nabla u\|_{L^2(\Om)}.
\eeno
Consequently, for $\epsilon>0$ to be determined later, \ben\label{F}
\|\sqrt{\rho}F\|^2_{L^2(\Om)}\leq \epsilon\|\nabla^2
w\|^2_{L^2(\Om)}+C_{\epsilon}\big(1+\|\nabla w\|^2_{L^2(\Om)}
+\|\nabla u\|^2_{L^2(\Om)}\big). \een
Noting that $Lw=\rho\partial_t w-\rho F$, we get by using Proposition \ref{proL} again that
\beno
\|\nabla^2w\|^2_{L^2(\Om)}\leq C\big(\|\rho\partial_t
w\|^2_{L^2(\Om)}+\|\rho F\|^2_{L^2(\Om)}\big)\leq
C\big(\|\sqrt{\rho}\partial_t w\|^2_{L^2(\Om)}+\|\sqrt{\rho}
F\|^2_{L^2(\Om)}\big), \eeno
which implies by taking $\epsilon=\f{1}{3C}$ in (\ref{F}) that
\beno
\|\sqrt{\rho}F\|^2_{L^2(\Om)}\leq \f{1}{2}\|\sqrt{\rho}\partial_t
w\|^2_{L^2(\Om)} +C\big(1+\|\nabla w\|^2_{L^2(\Om)}+\|\nabla
u\|^2_{L^2(\Om)}\big). \eeno

Substituting  this estimate into (\ref{eq:w-energy}) and noting that $\|\nabla u(t)\|^2_{L^2(\Om)}\in
L^1(0,T)$, the estimate (\ref{estw1}) follows from Gronwall's
inequality.\ef

\bthm{Corollary}\label{Cor-uest} Under the assumption (\ref{M}), we
have
\beno \|\nabla
u\|_{L^{\infty}(0,T;L^2(\Om))},\,\,\|u\|_{L^{\infty}(0,T;L^6(\Om))},\,\,
\|\nabla u\|_{L^{2}(0,T;L^q(\Om))}\leq C,
\eeno
for any $q\in [2,6]$.
\ethm

\noindent{\bf Proof.}\,\,This can be deduced from Proposition
\ref{lemmaw1}, (\ref{estv}) and Sobolev embedding theorem.\ef

\section{High order a priori estimates for the effective viscous flux}

In this section, we will give high order regularity estimates for $w$.
This is possible if the initial data $(\rho_0,u_0)$ satisfies the
compatibility condition (\ref{CompC}). We still assume that $(\rho,u)$ is a strong solution of
(\ref{NS}) in $[0,T)$ and satisfies (\ref{M}).The energy estimates in this
section are motivated by the calculations of D. Hoff \cite{Hoff}.

We begin by introducing some notations.
For a function or vector field(or even a $3\times 3$ matrix) $f(t,x)$,
the material derivative $\dot{f}$ is defined by
$$\dot{f}\eqdefa f_t+u\cdot\nabla f,$$
and $\div(f\otimes u)\eqdefa \sum_{j=1}^3\partial_j(fu_j)$.
For two matrices $A=(a_{ij})_{3\times 3}$ and $B=(b_{ij})_{3\times 3}$, we use the notation $A:B=\sum_{i,j=1}^3a_{ij}b_{ij}$ and $AB$ is as usual the multiplication of matrix.

We rewrite the second equation of (\ref{NS}) as
\beno
\rho\dot{u}+\nabla p(\rho)-Lu=0.
\eeno
By taking the material derivative to the above equation and using the fact $\dot{f}=f_t+\div(fu)-f\div u$, we obtain
\begin{eqnarray}\label{eq41}
&&\rho \dot{u}_t+\rho u\cdot \nabla \dot{u}+\nabla p_{t}+\div(\nabla p\otimes u)\nonumber\\
&&\quad=\mu\big[\Delta u_t+\div(\Delta u\otimes u)\big]+
(\lambda+\mu)\big[\nabla\div u_{t}+\div((\nabla\div u)\otimes u)\big].
\end{eqnarray}

Multiplying (\ref{eq41}) by $\dot{u}$ and integrating on $\Omega$ to obtain
\ben\label{eq42}
&&\f{d}{dt}\int\limits_\Omega\f{1}{2}\rho|\dot{u}|^2dx-\mu\int\limits_{\Om}\dot{u}\cdot\big(\Delta
u_t+\div(\Delta u\otimes u)\big)dx\non\\
&&\quad-(\lambda+\mu)\int\limits_{\Om}\dot{u}\cdot\big((\nabla
\div u_{t})+\div((\nabla\div u)\otimes u))\big)dx\non\\
&&=\int\limits_\Omega p_{t}\div\dot{u}+(\dot{u}\cdot\nabla u)\cdot\nabla p dx.
\een
The $\mu$-term can be calculated as follows.
\beno
&&-\int\limits_{\Om}\dot{u}\cdot\big(\Delta
u_t+\div(\Delta u\otimes u)\big)dx=\int\limits_{\Om}\left[\nabla\dot{u}:\nabla u_t+ u\otimes\Delta u:\nabla \dot{u}\right]dx\\
&&=\int\limits_{\Om}\Big[|\nabla\dot{u}|^2-\nabla(u\cdot\nabla u):\nabla\dot{u}+ u\otimes\Delta u:\nabla \dot{u}\Big]dx\\
&&=\int\limits_{\Om}\Big[|\nabla\dot{u}|^2-\big((\nabla u\nabla u)+(u\cdot\nabla) \nabla u\big):\nabla\dot{u}-\nabla(u\cdot\nabla\dot{u}):\nabla u\Big]dx\\
&&=\int\limits_{\Om}\Big[|\nabla\dot{u}|^2-(\nabla u\nabla u):\nabla\dot{u}-\div(\nabla u\otimes u):\nabla\dot{u}\\
&&\quad-(\nabla u\nabla\dot{u}):\nabla u-((u\cdot\nabla)\nabla\dot{u}):\nabla u\Big]dx\\
&&=\int\limits_{\Om}\Big[|\nabla\dot{u}|^2-(\nabla u\nabla u):\nabla\dot{u}+((u\cdot\nabla)\nabla\dot{u}):\nabla u\\
&&\quad-(\nabla u\nabla\dot{u}):\nabla u-((u\cdot\nabla)\nabla\dot{u}):\nabla u\Big]dx\\
&&\geq \int\limits_{\Om}\left[\f{3}{4}|\nabla\dot{u}|^2-C|\nabla u|^4\right]dx.
\eeno

To estimate the $(\lambda+\mu)$-term of (\ref{eq42}), note that
\beno
&&\div((\nabla\div u)\otimes u)=\nabla(u\cdot\nabla\div u)-\div(\div u\nabla\otimes u)+\nabla(\div u)^2,\\
&&\div\dot{u}=\div u_t+\div(u\cdot\nabla u)=\div u_t+u\cdot\nabla\div u+\nabla u:(\nabla u)'.
\eeno
Here $A'$ means the transpose of matrix $A$. We have
\beno
&&-\int\limits_{\Om}\dot{u}\cdot\Big[\nabla\div
u_{t}+\div((\nabla\div u)\otimes u)\Big]dx\\
&=&\int\limits_\Omega\Big[\div\dot{u}\div u_t+\div\dot{u}(u\cdot\nabla\div u)\\
&&\quad-\div u(\nabla\dot{u})':\nabla u+\div\dot{u}(\div u)^2\Big]dx\\
&=&\int\limits_\Omega\Big[|\div\dot{u}|^2-\div\dot{u}\nabla u:(\nabla u)'-\div u(\nabla\dot{u})':\nabla u+\div\dot{u}(\div u)^2\Big]dx\\
&\geq&\int\limits_\Omega\Big[\f{1}{2}|\div\dot{u}|^2-\f{1}{4} |\nabla\dot{u}|^2-C|\nabla u|^4\Big]dx.
\eeno

We continue to estimate the pressure term.
\beno
&&\int\limits_\Omega p_{t}\div\dot{u}+(u\cdot\nabla\dot{u})\cdot\nabla p dx
\\&&=\int\limits_\Omega p'(\rho)\rho_t \div\dot{u}+({u}\cdot\nabla\dot u)\cdot\nabla p dx\\
&&=\int\limits_\Omega -\rho p'(\rho)\div u\div\dot{u}-(u\cdot\nabla p(\rho))\div\dot{u}+({u}\cdot\nabla\dot u)\cdot\nabla p dx\\
&&=\int\limits_\Omega -\rho p'(\rho)\div u\div\dot{u}+p\Big[\div((\div\dot{u})u)-\div(({u}\cdot\nabla\dot u))\Big]dx\\
&&=\int\limits_\Omega -\rho p'(\rho)\div u\div\dot{u}+p\Big[\div u\div\dot{u}-(\nabla u)':\nabla\dot{u}\Big]dx\\
&&\leq C\|\nabla u\|_{L^2(\Om)}\|\nabla\dot{u}\|_{L^2(\Om)}\leq C\|\nabla\dot{u}\|_{L^2(\Om)},
\eeno
where we used the assumption (\ref{M}) and Corollary \ref{Cor-uest} in the last two inequalities.

Substituting those estimates into (\ref{eq42}) yields
\ben\label{eq43}
&&\f{d}{dt}\int\limits_\Omega\rho|\dot{u}|^2dx+\mu\int\limits_{\Om}|\nabla\dot{u}|^2dx
+(\lambda+\mu)\int\limits_{\Om}|\div\dot{u}|^2dx\non\\
&&\leq C\int\limits_{\Om}|\nabla{u}|^4dx+C\|\nabla\dot{u}\|_{L^2(\Om)}.
\een

To conclude the estimate by Gronwall's inequality, we will use the term $\|\sqrt{\rho}\dot{u}\|_{L^2(\Om)}$ to control $\|\nabla u\|_{L^4(\Om)}$.
Thanks to the definition of $w$, we know that $w$ satisfies
\ben\label{eqforw}
\mu\Delta w+(\lambda+\mu)\nabla \div w=\rho \dot{u}\hspace{4pt}\text{in}\hspace{4pt} \Omega,
\een
with the zero boundary condition. We get by Proposition 2.1 that
\begin{eqnarray*}
\|\nabla^2 w\|_{L^2(\Omega)}\leq C\|\rho \dot{u}\|_{L^2(\Omega)}\leq C\|\sqrt{\rho}\dot{u}\|_{L^2(\Omega)},
\end{eqnarray*}
which together with the interpolation inequality, Corollary \ref{Cor-uest}, and Proposition 2.1 leads to
\begin{eqnarray*}
\|\nabla u\|^4_{L^4(\Omega)}&\leq&\|\nabla u\|_{L^2(\Omega)}\|\nabla
u\|^3_{L^6(\Omega)}\leq C\|\nabla u\|_{L^6(\Omega)}\|\nabla
u\|^2_{L^6(\Omega)}\\
&\leq& C\|\nabla u||^2_{L^6(\Omega)}\big(\|\nabla
w\|_{L^6(\Omega)}+\|\nabla v\|_{L^6(\Omega)}\big)\\
&\leq& C\|\nabla u||^2_{L^6(\Omega)}\big(1+\|\nabla^2
w\|_{L^2(\Omega)}\big)\\
&\leq& C\|\nabla u\|^2_{L^6(\Omega)}\Big(1+\|\sqrt{\rho}\dot{u}\|_{L^2(\Omega)}\Big).
\end{eqnarray*}

Substituting this estimate into (\ref{eq43})and noting that  $||\nabla
u(t)||^2_{L^6(\Omega)}\in L^1(0,T)$ by Corollary 3.3, we get by
Gronwall's inequality that
\begin{eqnarray}\label{eq:w-high}
\int\limits_\Omega\rho|\dot{u}|^2dx+\int^T_0\int\limits_\Omega|\nabla
\dot{u}|^2dxdt\leq C,
\end{eqnarray}
with $C$ depending only on $T,M$ and $\rho_0,u_0,g$. Here we used the compatibility condition (\ref{CompC}).
\vspace{0.1cm}

With the help of Sobolev embedding theorem and using  the equation (\ref{eqforw}) again,
we deduce from (\ref{eq:w-high}) that
\begin{proposition}\label{Prop41}
Under the assumption (3.2), we have for all $2\leq q\leq 6$,
\begin{equation}\label{estw2}
\|\nabla w\|_{L^2(0,T;L^{\infty}(\Omega))},\,\,\|\nabla^2w\|_{L^2(0,T;L^q(\Omega))}\leq C,
\end{equation}
with the constant C depending on $q,M, T$ and $\rho_0,
u_0, g $.
\end{proposition}

\section{Proof of Theorem \ref{main}}

Now we are in position to prove Theorem \ref{main}.
We will prove it by the contradiction argument. Assume that $T^*<\infty$ and
$$\sup\limits_{s\in [0,T^*)}\|\rho(s)\|_{L^{\infty}(\Om)}<\infty.$$
By Theorem \ref{thmCCSS}, it suffices to show that
\ben\label{eq:blow}
\sup\limits_{s\in [0,T^*)}\|\na \rho(s)\|_{L^q(\Om)}<\infty.
\een

Taking the derivative with respect to $x$ for the first equation of (\ref{NS}) to obtain
\ben\label{eq52}
\partial_t \nabla\rho+(u\cdot\nabla)\nabla \rho+\nabla u\nabla\rho+\div u\nabla\rho+\rho\nabla\div u=0.
\een
In the following estimates we will use
\ben
\|\nabla^2v\|_{L^q(\Omega)}&\leq& C\|\nabla\rho\|_{L^q(\Omega)},\label{eq:v1}\\
\|\nabla v\|_{L^\infty(\Omega)}&\leq&
C\Big(1+\|\nabla v\|_{BMO(\Om)}\ln(e+\|\nabla^2v\|_{L^q(\Omega)})\Big)\non\\
&\leq& C\Big(1+\|\rho\|_{L^{\infty}\cap L^2(\Om)}\ln(e+\|\nabla\rho\|_{L^q(\Omega)})\Big)\non\\
&\leq& C\Big(1+\ln(e+\|\nabla\rho\|_{L^q(\Omega)})\Big)\label{eq:v2}
\een
with the second estimate followed from Proposition \ref{proL}, \ref{propend} and Lemma \ref{embedding2}.

Multiplying  (\ref{eq52}) by $q|\nabla\rho|^{q-2}\nabla\rho$
and integrating the resulting equation on $\Om$, we obtain
\beno
\frac{d}{dt}\int\limits_{\Om}|\nabla\rho|^qdx&\leq& C\int\limits_{\Om}|\nabla u||\nabla\rho|^qdx+q\int\limits_{\Om}\rho|\nabla\div u||\nabla\rho|^{q-1}dx\non\\
&\leq& C\|\nabla
u\|_{L^{\infty}(\Om)}\|\nabla\rho\|^q_{L^q(\Om)}+C\|\nabla^2u\|_{L^q(\Om)}\|\nabla\rho\|^{q-1}_{L^q(\Om)}\non\\
&\leq& C\big(\|\nabla w\|_{L^{\infty}(\Om)}+\|\nabla
v\|_{L^{\infty}(\Om)}\big)\|\nabla\rho\|^q_{L^q(\Omega)}\\&&\quad+C\big(\|\nabla^2w\|_{L^q(\Omega)}
+\|\nabla^2v\|_{L^q(\Omega)}\big)\|\nabla\rho\|^{q-1}_{L^q(\Omega)},
\eeno
from which and (\ref{eq:v1})-(\ref{eq:v2}), we infer that
\beno
\frac{d}{dt}\int\limits_{\Om}|\nabla\rho|^qdx
&\le& C\big(1+\|\nabla v\|_{L^{\infty}(\Om)}+\|\nabla w\|_{L^{\infty}(\Om)}\big)\|\nabla\rho\|^q_{L^q(\Omega)}\\
&&\quad+C\|\nabla^2w\|_{L^q(\Omega)}\|\nabla\rho\|^{q-1}_{L^q(\Omega)}\non\\
&\leq& C\big(1+\|\nabla w\|_{W^{1,q}(\Omega)}+\ln(e+\|\nabla\rho\|_{L^q(\Omega)})\big)\|\nabla\rho\|^q_{L^q(\Omega)}\non\\
&&\quad+\|\nabla^2w\|_{L^q(\Omega)}\|\nabla\rho\|^{q-1}_{L^q(\Omega)}.\non
\eeno
Note that $\|\nabla w\|_{W^{1,q}(\Omega)}\in L^2(0,T^*)$ by Proposition \ref{Prop41}.
Then by Gronwall's inequality, we conclude the proof of (\ref{eq:blow}) and hence Theorem \ref{main}.\ef

\section*{Acknowledgement}
The authors thank Professors Song Jiang, Changxing Miao, Zhouping Xin and Ping Zhang
for their profitable discussion and suggestions. This work was done while Yongzhong Sun and Zhifei Zhang were
visiting the Morningside Center of Mathematics in CAS. We would like
to thank the hospitality and support of the Center. Yongzhong Sun is
supported by NSF of China under Grant 10771097. Zhifei Zhang is
supported by NSF of China under Grant 10990013.

\end{document}